\documentclass[11pt]{article}
 \usepackage{amssymb,amsmath,amsthm}

\newtheorem{df}{Definition}[section]
\newtheorem{prop}[df]{Proposition}
\newtheorem{lemma}[df]{Lemma}
\newtheorem{teo}[df]{Theorem}

\theoremstyle{definition}

\newtheorem{ex}[df]{Example}

\newcommand{\sezione}[1]{\section{#1}\setcounter{equation}{0}}

\def\R{\mathbb{R}}
\def\N{\mathbb{N}}

\def\e{{\varepsilon}}
\def\eps{{\varepsilon}}

\def\span{\mathop{\rm span}\nolimits}

\def\sotto{\mathop{ \underbrace {\pm\gamma_i,\dots,\pm\gamma_i}  }\nolimits}

\def\X{{\rm X}}
\def\to{\rightarrow}
\def\la{\lambda}

 \newcommand{\uni}{\mathop{\cup}}
\newcommand{\inte}{\mathop{\cap}}

\begin{document}


\title{ \huge\textsc \bf On the exact number of bifurcation branches
from a multiple eigenvalue \footnote{The first author is supported
by the M.I.U.R. National Project ``Metodi variazionali ed equazioni
differenziali nonlineari''. The second author is supported by the
M.I.U.R. National Project ``Metodi variazionali e topologici nello
studio di fenomeni non lineari''. } }
\author{\small\sc{
Dimitri MUGNAI \footnote{Dipartimento di Matematica e Informatica,
Universit\`a di Perugia, via Vanvitelli 1, 06123 Perugia, Italy,
e-mail: mugnai@dipmat.unipg.it.},\ Angela PISTOIA
\footnote{Dipartimento di Metodi e Modelli Matematici,
Universit\`{a} di Roma "La Sapienza", 00100 Roma, Italy, e-mail:
pistoia@dmmm.uniroma1.it.} } }
 \date{\today }

  \maketitle


\begin{abstract}
We study local bifurcation from an eigenvalue with multiplicity
greater than one for a class of semilinear elliptic equations. We
evaluate the exact number of bifurcation branches of non trivial
solutions and we compute the Morse index of the solutions in those
branches.

\end{abstract}


{\small
\noindent{\bf Key words}:
local bifurcation, multiple branches, multiple eigenvalue, Morse index.


\noindent{\bf A. M. S. subject classification 2000}:
 35B32, 35J20, 35J60.

}


\sezione{Introduction and main results}

Let us consider the problem
\begin{equation}
\label{p1}
  \begin{cases}
  -\Delta u=|u|^{p-1}u+\lambda u &  \mbox{ in }\Omega\\
u=0 &  \mbox{ on }\partial\Omega,
\end{cases}
\end{equation}
where $\Omega$ is a bounded open domain in $\R^N$, $p>1$ and
$\lambda\in\R.$

It has the trivial family of solutions $\{(\lambda,0)\ |\
\lambda\in\R\}.$   A point $(\lambda^*,0)$ is called a {\it
bifurcation point} for  (\ref{p1})  if every neighborhood of
$(\lambda^*,0)$ contains nontrivial solutions of (\ref{p1}). It is
easily seen that a necessary condition  for $(\lambda^*,0)$  to be
a bifurcation point is that $\lambda^*$ is an eigenvalue of the
problem
\begin{equation}\label{eip}
\begin{cases}-\Delta u=\lambda u& \hbox{in}\ \Omega,\\ u=0& \hbox{on}\
\partial\Omega.\end{cases}
\end{equation}
We denote by $\lambda_1<\lambda_2\le\dots\le \lambda_j\le\dots $
the sequence of the eigenvalues of the problem (\ref{eip}).

Since problem (\ref{p1}) has a variational structure, the fact that
$\lambda^*$ is an eigenvalue of the problem (\ref{eip}) is not only
necessary, but  is also a sufficient condition for bifurcation to
occur. More precisely in \cite{B} and \cite{M}, it was proved that
for any eigenvalue $\lambda_j$ of (\ref{eip}) there exists $r_0>0$
such that for any $r\in(0,r_0)$ there are at least two distinct
solutions $(\lambda_i(r), u_i(r)),$ $i=1,2$ of (\ref{p1}) having
$\|u_i(r)\|=r$ and in addition $(\lambda_i(r),u_i(r)) \rightarrow
(\lambda_j,0)$ as $r\rightarrow0$.

A very interesting problem is to find the structure of the
bifurcation set at any eigenvalue $\lambda_j,$ namely the set of
nontrivial solutions $(\lambda,u)$ of (\ref{p1}) in a neighborhood
of $\lambda_j.$ In \cite{CR} the authors provide an accurate
description in the case of a simple eigenvalue, by showing that the
bifurcation set is a $C^1$ curve crossing $(\lambda_j,0)$. If the
eigenvalue $\lambda_j$ has higher multiplicity, in \cite{R} the
author describes the possible behavior of the bifurcating set by
showing that the following alternative occurs: either
$(\lambda_j,0)$ is not an isolated solution of (\ref{p1}) in
$\{\lambda_j\}\times H^1_0(\Omega)$, or there is a one-sided
neighborhood $U$ of $ \lambda_j$ such that for all $\lambda\in
U\setminus\{\lambda_j\}$ problem (\ref{p1}) has at least two
distinct nontrivial solutions, or there is a neighborhood $I$ of $
\lambda_j$ such that for all $\lambda\in I\setminus\{\lambda_j\}$
problem (\ref{p1}) has at least one nontrivial solution. We also
quote the results obtained in \cite{A}, where the author gives a
more precise description of the bifurcation set.

A natural question arises: {\it which is the exact number of
nontrivial solutions of \eqref{p1} bifurcating from any eigenvalue
$\lambda_j?$}

This question is the object of the present paper.

Let us fix an eigenvalue $\lambda_j$ with multiplicity $k$, i.e.
$\lambda_{j-1}<\lambda_j=\dots=\lambda_{j+k-1}<\lambda_{j+k}\le
\dots $. In view of the discussion above, we are mainly interested
in the case $j\geq 2$ and $k\geq 2$, though our results hold also
for $j=1$ or $k=1$.

We state our first main result.

\begin{teo}\label{main1}
There exists $\delta=\delta(\la_j)>0$ such that for any
$\lambda\in(\lambda_j-\delta,\lambda_j)$ problem \eqref{p1} has at
least $k$ pairs of solutions $(u_\lambda,-u_\lambda)$ bifurcating
from $\lambda_j$. Moreover, associated to each $u_\lambda$ there exist
real numbers $a^1_\lambda,\dots,a^k_\lambda$ and a function
$\phi_\lambda\in {\rm X}$ (see \eqref{ics1} and \eqref{ics2}),
with $(\phi_\lambda,e_i)=0$ for $i=1,\dots,k$ such that
\begin{equation}\label{ul}
u_\lambda
=(\lambda_j-\lambda)^{\frac{1}{p-1}}\left[\sum\limits_{i=1}^ka^i_\lambda
e_i +\phi_\lambda\right],\ \lim\limits_{\lambda\rightarrow
\lambda_j}\| \phi_\lambda\|_{\rm X}=0\ \hbox{and}\
\lim\limits_{\lambda\rightarrow\lambda_j}a^i_\lambda=a^i,
\end{equation}
where $a:=(a^1,\dots, a^k)$ is a critical point of the function
$J_{\lambda_j}:\R^k\longrightarrow\R$, defined by
\begin{equation}\label{J}
J_{\lambda_j}(a)=\frac{1}{2}\left(a_1^2+\dots+a_k^2\right)-
\frac{1}{p+1}\int\limits_\Omega\left|a_1e_1
+\dots+a_ke_k\right|^{p+1}dx.
\end{equation}
Here $e_1,\dots,e_k,$ are $k$ orthogonal eigenfunctions
associated to the eigenvalue $\lambda_j $ such that
$\|e_i\|_{L^2}=1$ for any $i=1,\dots,k$.
\end{teo}

We point out that the existence of at least $2k$  nontrivial
solutions bifurcating from the eigenvalue $\lambda_j$ was already
known (see \cite{B}, \cite{M}, \cite{R} and \cite{CFS}, \cite{GR}
and the references therein for the critical case). But in Theorem
\ref{main1} we also describe the asymptotic behaviour of those
solutions as $\lambda$ goes to $\lambda_j.$ In particular we find
out a relation between solutions  to problem \eqref{p1}
bifurcating from the eigenvalue $\lambda_j$ and critical points of
the function $J_{\lambda_j}$: the solution $u_\lambda$ which
satisfies \eqref{ul} is ``generated" by the critical point $a$.
This suggests that the solution  $u_\lambda$ ``generated" by  $a$
can inherit some properties of $a.$ More precisely we prove the
following result.

\begin{teo}\label{main2}
Assume $a$ is a non degenerate critical point of the function
$J_{\lambda_j}$ with Morse index $m.$ Then there exists
$\delta=\delta(\la_j)>0$ such that for any
$\lambda\in(\lambda_j-\delta,\lambda_j)$ problem \eqref{p1} has a
unique solution $ u_\lambda $  bifurcating from $\lambda_j$ which
satisfies \eqref{ul}. Moreover $u_\lambda$ is non degenerate and
its Morse index is $m+j-1.$
\end{teo}

Moreover,  it also suggests that the number of nontrivial
solutions to \eqref{p1} bifurcating from the eigenvalue
$\lambda_j$ coincides with the number of nontrivial critical
points of $J_{\lambda_j}.$ In fact we have the following result.

\begin{teo}\label{main3}
Assume the function $J_{\lambda_j}$ defined in \eqref{J} has
exactly $2h$ non trivial critical points which are non degenerate.
Then there exists $\delta=\delta(\la_j)>0$ such that for any
$\lambda\in(\lambda_j-\delta,\lambda_j)$ problem \eqref{p1} has
exactly $h$ pairs of solutions $(u_\lambda,-u_\lambda)$
bifurcating from $\lambda_j$ and all of them are non degenerate.
\end{teo}

The last theorem allows us to  give an accurate description of the
bifurcation branches from any multiple eigenvalues when $\Omega$
is a rectangle in $\R^2$ (see Example \ref{rect}) and also from
the second eigenvalue when $\Omega$ is a cube in $\R^3$ (see
Example \ref{cubo}), provided $p=3.$ We point out that when
$\Omega=(0,\pi)\times(0,\pi)$ is a square in $\R^2$ and $p$ is an
odd integer, in \cite{DGM} the authors studied the bifurcation
from the second eigenvalue $\lambda_2$ and they proved that the
bifurcation set is constituted exactly by the union of four $C^1$
curves crossing $(\lambda_2,0)$ from the left.
We also remark that some exactness results in bifurcation theory for a different class of 
problems were obtained in  \cite{Z}.

The proof of Theorem \ref{main1} is based upon a Ljapunov-Schmidt
reduction method (see, for example, \cite {A}, \cite{Ba}, \cite{MP}, \cite{R}).
The proofs of Theorem \ref{main2} and Theorem \ref{main3} rely on
the asymptotic estimate \eqref{ul}, which links nontrivial solutions
of \eqref{p1} bifurcating from $\lambda_j$ and critical points of
the function $J_{\lambda_j}.$

The paper is organized as follows: in Section 2 we introduce some
notation, in Section 3 we reduce problem \eqref{p1} to a finite
dimensional one, in Section 4 we study the reduced problem and we
prove Theorem \ref{main1}, in Section 5 we prove Theorem
\ref{main2} and Theorem \ref{main3} and in Section 6 we discuss
some applications.

In order to make the reading more fluent, in many calculations we
have used the symbol $c$ to denote different absolute constants
which may vary from line to line.


\sezione{Setting of the problem}

First of all we   rewrite problem (\ref{p1}) in a different way.
We introduce a positive parameter $\eps.$ An easy computation
shows that, if $u(x)$ solves problem (\ref{p1}), then for any
$\e>0$ the function $v(x)=\e^{-\frac{1}{p-1}}u(x)$ solves
\begin{equation}
\label{p2}
\begin{cases}-\Delta v=\e|v|^{p-1}v+\lambda v &  \mbox{ in
 }\Omega \\
v=0 &  \mbox{ on }\partial\Omega.
\end{cases}
\end{equation}

The parameter $\eps$ will be chosen in Lemma \ref{jj} as
$\e=\lambda_j-\lambda>0.$

Let ${\rm H}^1_0( \Omega)$ be the Hilbert space  equipped with the
usual inner product $\langle u,v \rangle=\int\limits_{\Omega
}\nabla u\nabla v$, which induces the standard norm $\Vert u\Vert
=\big(\ \int\limits_{\Omega }|\nabla u|^2\big)^{1/2}.$

If $r\in [1,+\infty)$ and $u\in {\rm L}^r(\Omega)$, we will set
$\|u\|_r=\big(\int\limits_\Omega |u|^r\big)^{1/r}$.

\begin{df}
\label{istar} Let us consider the embeddings $i :{\rm
H}^1_0(\Omega )\hookrightarrow{\rm L}^{\frac{2N}{N-2}}(\Omega ) $
if $N\ge3$ and $i :{\rm H}^1_0(\Omega )\hookrightarrow
\inte\limits_{q>1}{\rm L}^{q}(\Omega ) $ if $N=2.$ Let $i^*:{\rm
L}^{\frac{2N}{N+2}}(\Omega )\longrightarrow{\rm H}^1_0(\Omega ) $
if $N\ge3$ and $i^* : \uni\limits_{q>1}{\rm L}^{q}(\Omega
)\longrightarrow {\rm H}^1_0(\Omega )$ if $N=2$ be the adjoint
operators defined by
$$i^*(u)=v\quad\Longleftrightarrow\quad
\langle v,\varphi\rangle=\int\limits_{\Omega }
u(x)\varphi(x)dx\quad \forall\ \varphi\in{\rm H}^1_0(\Omega) .$$
\end{df}

It holds
\begin{eqnarray}\label{istar2}
& & \Vert{i^*(u)}\Vert \le c\Vert
u\Vert_{\frac{2N}{N+2}}\ \hbox{for any}\ u\in{\rm
L}^{\frac{2N}{N+2}}(\Omega ),\ \hbox{if}\ N\ge3, \\
& & \Vert{i^*(u)}\Vert \le c(q)\Vert
u\Vert_{q }\ \hbox{for any}\ u\in{\rm
L}^{q}(\Omega ),\ q>1,\  \hbox{if}\ N=2.
\end{eqnarray}
Here the positive constants $c$ and $c(q)$ depend only on $\Omega$
and $N$ and $\Omega,$ $N$ and $q,$ respectively.

Let us recall the following regularity result proved in \cite{ADN},
which plays a crucial role when $p>\frac{N+2}{N-2}$ and $N\ge3$.

\begin{lemma}
\label{adn2} Let $N\ge3$ and $s>\frac{2N}{N-2}.$ If $u\in {\rm
L}^{\frac{Ns}{N+2s}}(\Omega) $, then $i^*(u)\in{\rm L}^{s}(\Omega)$
and $\Vert i^*(u)\Vert_{s}\le c\Vert u\Vert_{\frac{Ns}{N+2s}}$,
where the positive constant $c$ depends only on $\Omega$, $N$ and
$s.$
\end{lemma}

Now, we introduce the space
\begin{equation}\label{ics1}
{\rm X} = {\rm H}^1_0(\Omega)\  \hbox{if either $N=2$ or $N\ge3$
and $1<p\le \frac{N+2}{N-2}$}
\end{equation}
or
\begin{equation}\label{ics2}
{\rm X} = {\rm H}^1_0(\Omega)\cap{\rm L}^{s}(\Omega),\
s=\frac{N(p-1)}{2},\ \hbox{if $N\ge3$ and $p>\frac{N+2}{N-2}$}.
\end{equation}
We remark that the choice of $s$ is such that
$\frac{pNs}{N+2s}=s$, a fact that will be used in the following.

${\rm X} $ is a Banach space equipped with the norm $\Vert u\Vert_\X
=\Vert u\Vert  $ in the first case and $\Vert u\Vert_\X =\Vert
u\Vert +\Vert u\Vert_{s} $ in the second case.

By means of the definition of the operator $i^{*}$, problem
(\ref{p2}) turns out to be equivalent to
\begin{equation}\label{p3}
\begin{cases}
u=i^*[\eps f(u ) +\lambda u] \\
u\in {\rm X} ,\end{cases}
\end{equation}
where $f(s)=|s|^{p-1}s$.

Now, let us fix an eigenvalue $\lambda_j$ with multiplicity $k$,
i.e. $\lambda_{j-1}<\lambda_j=\dots=\lambda_{j+k-1}<\lambda_{j+k}\le
\dots $.  We denote by $e_1,\dots,e_k,$ $k$ orthogonal
eigenfunctions associated to the eigenvalue $\lambda_j $ such that
$\|e_i\|_2=1$ for $i=1,\dots,k.$ We will look for solutions to
(\ref{p2}), or to (\ref{p3}), having the form
\begin{equation}\label{sol}
u(x)=\sum\limits_{i=1}^ka^i_\lambda e_i(x)+ \phi_\lambda
(x)=a_\lambda e+\phi_\lambda,
\end{equation}
where $a^i_\lambda \in\R$, the function $\phi_\lambda$ is a lower
order term and we have set $a:=(a^1,\dots,a^k),$
$e:=(e_1,\dots,e_k)$ and $a e:=\sum\limits_{i=1}^ka^i e_i.$

We consider the subspace of ${\rm X}$ given by $K_j= \span\left\{
 e_i\ |\ i= 1,\ldots,k\right\}
$
and its complementary space
$
K_j^{\perp}=\left\{\phi\in {\rm X}\ |\
\langle \phi,e_i\rangle=0,\   i=1,\dots,k\right\}.
$

Moreover let us introduce the operators $\Pi_j:{\rm X}\rightarrow
K_j$ and $\Pi_j^\perp:{\rm X}\rightarrow K_j^\perp$ defined by $
\Pi_j (u)=\sum_{i=1}^k  \langle u,e_i\rangle e_i $ and $\Pi^\perp_j
(u)=u-\Pi_j (u).$ We remark that there exists a positive constant $
c$ such that
\begin{equation}\label{civuole}
\Vert{\Pi_j(u)}\Vert_\X\le c\Vert u\Vert_\X,\quad
\|\Pi_j^\perp(u)\|_{\rm X}\leq c\|u\|_{\rm X} \quad
\forall\,u\in{\rm X}.
\end{equation}

Our approach to solve problem (\ref{p3}) will be to find, for
$\lambda$ close enough to $\lambda_j $ and $\e$ small enough, real
numbers $a^1,\dots,a^k$ and a function $\phi\in K_j^{\perp}$ such
that
\begin{equation}
\label{1}
\Pi_j^{\perp}\left\{ae +\phi-
i^*\left [\e f\left(ae+\phi\right)+\lambda(ae+\phi)\right]\right\}=0
\end{equation}
and
\begin{equation}
\label{2} \Pi_j\left\{ae +\phi- i^*\left [\e
f\left(ae+\phi\right)+\lambda(ae+\phi)\right]\right\}=0.
\end{equation}


\sezione{Finite dimensional reduction}

In this section we will solve equation \eqref{1}. More precisely,
we will prove that for any $a\in\R^k,$ for   $\lambda$ close
enough to $\lambda_j$ and $\e$ small enough, there exists a unique
$\phi\in K_j^{\perp}$ such that (\ref{1}) is fulfilled.

Let us introduce the linear operator $L_\lambda:K^{\perp}_j\to
K^{\perp}_j$ defined by $L_\lambda(\phi)= \phi-\Pi^\perp_j
\left\{i^*\left[\lambda\phi \right]\right\}.$

We can prove the following result.

\begin{lemma}
\label{lin} There exists $\delta>0$ and a constant $c>0$ such that
for any $\lambda\in(\lambda_j-\delta,\lambda_j+\delta)$, the
operator $L_\lambda$ is invertible and it holds
\begin{equation}
\label{lep}
\|L_\lambda(\phi)\|_\X\geq c\|\phi\|_\X\qquad \forall\ \phi\in
K_j^{\perp}.
\end{equation}
\end{lemma}
\begin{proof}
First of all, we remark that $L_\lambda$ is surjective.

Concerning the estimate, we prove our claim when $N\ge3$ and
$p>\frac{N+2}{N-2}$, and we argue in a similar way in the other
cases.

Assume by contradiction that there are sequences $\delta_n\to 0$,
$\lambda_n\to\lambda_j$ and $\phi_n\in K_j^\perp$ such that
\[
\|L_{\lambda_n}(\phi_n)\|_\X<\frac{1}{n}\|\phi_n\|_\X.
\]

Without loss of generality we can assume
\begin{equation}\label{fin1}
\|\phi_n\|_\X=1\qquad  \mbox{ for any $n\in \N$}.
\end{equation}

If $h_n:=L_{\la_n}(\phi_n)\in \Pi_j^\perp$, then
\begin{equation}\label{hn0}
\|h_n\|_\X\to 0
\end{equation}
and
\begin{equation}\label{hn}
\phi_n-i^\ast(\la_n\phi_n)=h_n-\Pi_j\{i^\ast[\la_n\phi_n]\}=h_n+w_n,
\end{equation}
where $w_n\in K_j$.

First of all we point out  that $w_n=0$ for any $n\in \N$. Indeed,
multiply equation \eqref{hn} by $e_i,\,i=1,\ldots,k$, so that $
\langle w_n,e_i\rangle=-\la_n\int\limits_\Omega \phi_ne_i=0, $ so
that $w_n\in K_j^\perp$, and then $w_n=0$.

By \eqref{fin1}, we can assume that, up to a subsequence,
$\phi_n\rightarrow \phi$ weakly in $\X$ and strongly in ${\rm
L}^q(\Omega)$ for any $q\in\big[1,\frac{2N}{N-2}\big)$.
Multiplying \eqref{hn} by a test function $v$, we get
\[
\langle \phi_n,v\rangle -\la_n\int\limits_\Omega \phi_n
v\,dx=\langle h_n,v\rangle,
\]
and passing to the limit, by \eqref{hn0} we deduce that $\phi \in
K_j$. Since $\phi \in K_j^\perp$, we conclude that $\phi=0$.

On the other hand, multiplying \eqref{hn} by $\phi_n$, we get
\[
\langle \phi_n,\phi_n\rangle-\la_n \int\limits_\Omega
\phi_n^2=\langle h_n,\phi_n\rangle,
\]
which implies $ \Vert \phi_n\Vert\rightarrow 0.$ Moreover by
(\ref{hn}),  Lemma \ref{adn2} and by interpolation (since
$1<\frac{Ns}{N+2s}<s$), we deduce that for some $\sigma\in(0,1)$
$$\Vert \phi_n\Vert_s\le c\left(\Vert h_n\Vert_s +
\Vert \phi_n\Vert_\frac{Ns}{N+2s}\right)\le c\left(\Vert
h_n\Vert_s +\Vert \phi_n\Vert^\sigma\right)
$$
(recall that $\frac{Ns}{N+2s}=\frac{s}{p}$) and  so $ \Vert
\phi_n\Vert_s\rightarrow 0.$ Finally a contradiction arises, since
$ \Vert \phi_n\Vert_{\rm X}=1.$
\end{proof}

Now we can  solve Equation (\ref{1}).

\begin{prop}
\label{i} For any compact set $W$ in $\R^k$ there exist
$\eps_0>0,$ $\delta>0$ and $R>0$ such that, for any $a\in W,$ for
any $\e\in(0,\e_0)$ and for any
$\lambda\in(\lambda_j-\delta,\lambda_j+\delta)$, there exists a
unique $\phi_\lambda(a)\in K_j^\perp$ such that
\begin{equation}
\label{i1} \Pi_j^{\perp}\left\{ae +\phi_\lambda(a)- i^*\left [\e
f\left(ae+\phi_\lambda(a)\right)+
\lambda(ae+\phi_\lambda(a))\right]\right\}=0.
\end{equation}
Moreover
\begin{eqnarray}
\label{sti2} \|\phi_\lambda(a)\|_{\rm X} \le R\eps.
\end{eqnarray}
Finally, the map $a\mapsto \phi_\lambda(a)$ is an odd
$C^1-$function from $\R^k$ to $K_j^\perp$.
\end{prop}
\begin{proof}
We prove our claim when $N\ge3$ and $p>\frac{N+2}{N-2}.$ We argue
in a similar way in the other cases.

Let us introduce the operator $T:K_j^\perp\longrightarrow
K_j^\perp$ defined by
$$
T(\phi):=\left(L_\lambda^{-1}\circ\Pi^\perp_j\circ i^*\right)
\left[\e f(ae+\phi)\right].
$$
We point out that $\phi$ solves equation (\ref{i1}) if and only
if $\phi$ is a fixed point of $T,$ i.e. $T(\phi)=\phi.$

Then, we will prove that there exist $\eps_0>0,$ $\delta>0$ and
$R>0$ such that, for any $a\in W,$ for any $\e\in(0,\e_0)$ and for
any $\lambda\in(\lambda_j-\delta,\lambda_j+\delta)$
$$
T:\{\phi\in K_j^\perp\ |\ \|\phi\|_\X\le
  R\eps\}\longrightarrow\{\phi\in K_j^\perp\ |\ \|\phi\|_\X\le
  R\eps\}
$$
is a contraction mapping.

First of all, let us point out that by Lemma \ref{lin},
\eqref{civuole}, (\ref{istar2}) and Lemma \ref{adn2}, we get that
there exists $c=c(N,s,\Omega,W)>0$ such that for any $\phi \in
K_j^\perp$, $a\in W$
\begin{eqnarray*}
& &\|T(\phi)\|_\X\leq c\eps \left[\|f(ae+\phi)\|_{\frac{2N}{N+2}}+
\|f(ae+\phi)\|_{\frac{Ns}{N+2s}} \right]\\
&\mbox{(H\"older inequality)} &\leq  c\eps
\|f(ae+\phi)\|_{\frac{Ns}{N+2s}} \leq c\eps
\left( \|ae\|^p_{\frac{Nsp}{N+2s}}+\|\phi\|^p_{\frac{Nsp}{N+2s}}\right)\\
& &\leq c\eps (1+\|\phi\|_\X^p).
\end{eqnarray*}

Finally, provided $\eps$ is small enough and $R$ is suitable chosen,
$T$ maps $\{\phi\in K_j^\perp\, :\, \|\phi\|_\X\le R\eps\}$ into
itself.

Now, let us  show that $T$ is a contraction, provided $\eps$ is
even smaller. As before, by Lemma \ref{lin}, \eqref{civuole},
\eqref{istar2} and Lemma \ref{adn2}, we get that there exists
$c>0$ such that for any $\phi_1,\,\phi_2\in K_j^\perp$, $a\in K$
 \begin{eqnarray*}
& &\|T(\phi_1)-T(\phi_2)\|_\X\leq c\eps  \left[
  \|f(ae+\phi_1)-f(ae+\phi_2)\|_{\frac{2N}{N+2}}\right.\\
& & \hskip5truecm
\left.+
\|f(ae+\phi_1)-f(ae+\phi_2)\|_{\frac{Ns}{N+2s}}\right]
\\ & &\leq c\eps \|f(ae+\phi_1)-f(ae+\phi_2)\|_{\frac{Ns}{N+2s}}\\
& & \le c\eps\left(\|\phi_1-\phi_2\|_{\frac{Nsp}{N+2s}}+ \|\phi_1
\|^{p-1}_{\frac{Nsp}{N+2s}}\|\phi_1-\phi_2\|_{\frac{Nsp}{N+2s}}
+\|\phi_1-\phi_2\|^p_{\frac{Nsp}{N+2s}}\right) .
\end{eqnarray*}
Indeed, by the mean value theorem, it follows that  there exists
$\vartheta\in(0,1)$ such that $
f(ae+\phi_1)-f(ae+\phi_2)=f'\left(ae+\phi_1+\vartheta
(\phi_2-\phi_1)\right)(\phi_1-\phi_2).$

Finally $\|T(\phi_1)-T(\phi_2)\|_\X\leq c\eps\|\phi_1-\phi_2\|_{\rm
X}$ if $\|\phi_1\|_\X,\|\phi_2\|_\X\le  R\eps$ and our claim
immediately follows.

The oddness of the mapping $a\mapsto\phi_\lambda(a)$ i.e.
$\phi_\lambda(a)=-\phi_\lambda(-a)$, is a straightforward
consequence of the uniqueness of solutions of problem \eqref{i1}.

The regularity of the mapping can be proved using  standard
arguments.
\end{proof}


\sezione{The reduced problem}

In this section we will solve equation \eqref{2}. More precisely,
we will prove that if $\lambda$ is close enough to $\lambda_j$,
there exists $a_\lambda\in\R^k$ such that equation \eqref{2} is
fulfilled.

Let $I_\lambda:{\rm H}^1_0(\Omega)\longrightarrow\R$ be defined by
\begin{equation}
\label{jl} I_\lambda(u):=\frac{1}{2}\int\limits_{\Omega}|\nabla
u|^2\,dx-\frac{\lambda}{2}\int\limits_{\Omega}
u^2\,dx-\frac{\eps}{p+1}\int\limits_{\Omega}|u|^{p+1}\,dx.
\end{equation}
It is well known that critical points of $I_\lambda$ are solutions
of problem (\ref{p2}). Let us consider the reduced functional
$J_\lambda:\R^k\longrightarrow\R$ defined by
\begin{equation}
\label{jtl}
J_\lambda(a):=I_\lambda\left(ae+\phi_\lambda(a)\right),
\end{equation}
where $\phi_\lambda(a)$ is the unique solution of \eqref{i1}.

\begin{lemma}
\label{jt} A function $u_\lambda:=ae+\phi_\lambda(a)$ is a
solution to \eqref{p2} if and only if $a$ is a critical point of
$J_\lambda$.
\end{lemma}
\begin{proof}
We only point out that
$$\frac{\partial J_\lambda}{\partial
a_i}(a) =J'_\lambda \left(ae+\phi_\lambda(a)\right)
\left(e_i+\frac{\partial\phi_\lambda}{\partial a_i}(a)\right)
 =J'_\lambda \left(ae+\phi_\lambda(a)\right)\left(e_i\right),
$$
since $\phi_\lambda(a)$ solves equation \eqref{i1} and
$\frac{\partial\phi_\lambda}{\partial a_i}(a)\in K_j^\perp$. Then
the claim easily follows.
\end{proof}

>From now on we assume
\begin{equation}\label{epsi}
\eps:=\lambda_j-\lambda>0.
\end{equation}
\begin{lemma}\label{jj}
It holds
\begin{equation}\label{45}
J_\lambda(a)=(\lambda_j-\lambda)\left[
J_{\lambda_j}(a)+\Phi_\lambda(a)\right],
\end{equation}
where $J_{\lambda_j}$ is defined in \eqref{J} and
$\Phi_\lambda:\R^k\longrightarrow\R$ is an even $C^1-$function such
that
 $\Phi_\lambda$ goes to zero $C^1-$uniformly on compact sets of $\R^k$ as
$\la\to \la_j$.
\end{lemma}

\begin{proof}
Set $\phi:=\phi_\lambda(a).$ It holds (using (\ref{epsi}))
\begin{eqnarray*}
J_\lambda(a)& &=\frac{1}{2}\int\limits_\Omega|\nabla
(ae+\phi)|^2\,dx- \frac{\eps}{p+1} \int\limits_\Omega| ae+\phi
|^{p+1}\,dx-\frac{\lambda}{2}\int\limits_\Omega( ae+\phi )^{2}\,dx
\nonumber\\ &
&=\frac{1}{2}(\lambda_j-\lambda)\left(a_1^2+\dots+a_k^2\right)
-\frac{\eps}{p+1}\int\limits_\Omega\left|a_1e_1+\dots+a_ke_k\right|^{p+1}
\,dx\nonumber\\ & & +\frac{1}{2}\int\limits_\Omega|\nabla\phi
|^2\,dx -\frac{\lambda}{2}\int\limits_\Omega
\phi^{2}\,dx-\frac{\eps}{p+1}
\int\limits_\Omega\left[\left|ae+\phi\right|^{p+1}-
\left|ae \right|^{p+1}\right]\,dx\nonumber\\
& &=(\lambda_j-\lambda) \left[J_{\lambda_j}
 (a)+\Phi_\lambda(a)\right],
\end{eqnarray*}
where $J_{\lambda_j}$ is defined in (\ref{J}) and
\begin{eqnarray*}
 & & \Phi_\lambda(a):=
\frac{1}{\lambda_j-\lambda}\left[\frac{1}{2}\int\limits_\Omega|\nabla
  \phi |^2\,dx
-\frac{\lambda}{2}\int\limits_\Omega \phi^{2}\,dx\right]  \nonumber\\
& &\hskip2truecm-\frac{1}{p+1}
\int\limits_\Omega\left[\left|ae+\phi\right|^{p+1}
-\left|ae\right|^{p+1}\right]\,dx.
\end{eqnarray*}
Here we used the fact that $\phi_{\lambda_j}(a)=0$ $\forall\, a$,
as it is clear from \eqref{sti2}.

Of course $\Phi_\lambda$ is even and of class $C^1$. It remains to
prove that it goes to zero $C^1$--uniformly on every compact
subset $W$ of $\R^k$ as $\la\to \la_j$, that is, as $\eps\to 0$.

Let us fix a compact set $W$ in $\R^k.$ It is easy to check that
$$|\Phi_\la(a)|\leq \frac{c}{\eps}\|\phi\|^2+c\|\phi\|_{\rm X}\le c\eps,
\ \hbox{for any $a\in W$}.$$ Indeed, by the mean value theorem we
deduce that there exists $\vartheta\in(0,1)$ such that
$$
\frac{1}{p+1} \int\limits_\Omega\left[\left|ae+\phi\right|^{p+1}
-\left|ae\right|^{p+1}\right]\,dx=\int\limits_\Omega
f\left(ae+\vartheta\phi\right)\phi\,dx.
$$
Therefore $\Phi_\lambda$ goes to zero uniformly on $W$ as $\la\to
\la_j$, since $\boldsymbol{\|\phi\|_{\rm X}\le R\eps}$ by
\eqref{sti2}.

Now, let us prove that also $\nabla \Phi_\la$ goes to zero  as
$\la\to \la_j$ uniformly on $W$. Indeed, fix $i=1,\ldots,k$ and
evaluate
\begin{equation}\label{amezzo}
\begin{aligned}
\frac{\partial \Phi_\la(a)}{\partial a_i}=&
\frac{1}{\lambda_j-\lambda}\bigg[\int\limits_\Omega \nabla \phi
\cdot \frac{\partial \nabla \phi}{\partial a_i}\,dx
-\lambda\int\limits_\Omega \phi\frac{\partial \phi}{\partial
a_i}\,dx \bigg]\\ &
 -\int\limits_\Omega\left[|ae+\phi|^{p-1}(ae+\phi)
 \Big(e_i+\frac{\partial \phi}{\partial a_i}\Big)
-|ae|^{p-1}aee_i\right]\,dx.
\end{aligned}\end{equation}

By \eqref{i1}, for every $z\in K_j^\perp$ we have
\[
\int\limits_\Omega \nabla \phi\cdot \nabla
z\,dx-\lambda\int\limits_\Omega \phi z\,dx-\eps \int\limits_\Omega
|ae+\phi|^{p-1}(ae+\phi)z\,dx=0 .
\]
Then, taking  $z=\frac{\partial \phi}{\partial a_i}\in K_j^\perp,$
by (\ref{amezzo}) we deduce
$$
\frac{\partial \Phi_\la(a)}{\partial
a_i}=\int\limits_\Omega\left[f(ae+\phi)-f(ae) \right]e_i\,dx $$
and so
$$
\left|\frac{\partial \Phi_\la(a)}{\partial a_i}(a)\right|\leq
c\|\phi\|_{\rm X}\le c\eps, \ \hbox{for any $a\in W$}.
$$

Indeed, again by the mean value theorem we deduce that there exists
$\vartheta\in(0,1)$ such that
$$
\int\limits_\Omega\left[f(ae+\phi)-f(ae) \right]e_i\,dx=
\int\limits_\Omega f'\left(ae+\vartheta\phi\right) \phi e_i\,dx.$$

Therefore, also $\nabla \Phi_\lambda$ goes to zero uniformly on
$W$ as $\la\to \la_j$.
\end{proof}

\begin{prop}\label{jrab}
There exists $\delta>0$ such that for any
$\lambda\in\left(\lambda_j-\delta,\lambda_j\right)$ the function
$J_\lambda$ has at least $k$ pairs  $(a_\la,-a_\la)$  of distinct
critical points. Moreover $a_\lambda\rightarrow a$ as $\lambda$ goes
to $\lambda_j$ and $a$ is a critical point of $J_{\lambda_j}$ (see
\eqref{J}).
\end{prop}
\begin{proof} First of all, we note that $J_\lambda(0)=0$ and
also that $J_\lambda$ is an even function. Moreover
(see \eqref{J}), it is clear that there exist $R>r>0$ such that
$$
\inf\limits_{|a|=r}J_{\lambda_j}(a) >J_{\lambda_j}(0)=0>
\sup\limits_{|a|=R}J_{\lambda_j}(a).
$$
Therefore, by Lemma \ref{jj} we deduce that, if $\lambda$ is close
enough to $\lambda_j$, it holds
$$
\inf\limits_{|a|=r}J_\lambda(a) >
J_\lambda(0)=0>\sup\limits_{|a|=R}J_\lambda(a).
$$
Then $J_{\lambda_j}$ has at least $k$ pairs of distinct critical
points $(a_\la,-a_\la)$ in $B(0,R)$.  We can assume that $a_\la
\to a\in \overline{B(0,R)}$ as $\la \to \la_j$. By \eqref{45} we
get $ \nabla J_{\lambda_j}(a_\la)=\frac{1}{\e}\nabla
J_\la(a_\la)-\nabla \Phi_\la(a_\la)=-\nabla \Phi_\la(a_\la), $ and
since $\Phi_\la $ goes to zero $C^1-$ uniformly on
$\overline{B(0,R)}$ as $\la \to \la_j$, we get $ \nabla
J_{\lambda_j}(a)=0. $ That proves our claim.
\end{proof}

\begin{proof}[Proof of Theorem \ref{main1}]

The claim follows by Lemma \ref{jt} and  Proposition \ref{jrab}.
\end{proof}

\sezione{Some uniqueness results}

First of all we describe the asymptotic behaviour of the solution
$u_\lambda$ of problem \eqref{p1} bifurcating from the eigenvalue
$\lambda_j$ as $\lambda$ goes to $\lambda_j.$

\begin{prop}\label{nece}
Let $u_\lambda\in{\rm X}$ be a solution to problem \eqref{p1} such
that $\|u_\lambda\|_{\rm X}$ goes to zero as $\lambda$ goes to
$\lambda_j$. Then for any $\lambda$ sufficiently close to
$\lambda_j$ there exist $a_\lambda^1,\dots, a_\lambda^k\in\R$ and
$\phi_\lambda\in K_j^\perp$ such that
$$
u_\lambda=(\lambda_j-\lambda)^\frac{1}{p-1}\left(
\sum\limits_{i=1}^ka_\lambda^ie_i+\phi_\lambda\right),
$$
where $\phi_\lambda \rightarrow 0$ in X and
$a_\lambda^i\rightarrow a^i,$ $i=1,\dots,k$, as $\lambda$ goes to
$\lambda_j$. Moreover $a=(a^1,\dots,a^k)$ is a critical point of
$J_{\lambda_j}.$
\end{prop}

\begin{proof} We prove our claim  when $N\ge 3$ and $p>\frac{N+2}{N-2}.$
We argue in a similar way in the other cases.

The function $v_\lambda=\frac{u_\lambda}{\|u_\lambda\|_{\rm X}}$
solves the problem
\begin{equation}\label{e1}
\begin{cases}
-\Delta v= \lambda v +\|u_\lambda\|_{\rm X}^{p-1}|v|^{p-1}v &
\mbox{ in }\Omega\\ v=0 &  \mbox{ on }\partial\Omega.
\end{cases}
\end{equation}
Up to a subsequence, we can assume that, as $\lambda$ goes to
$\lambda_j$, $v_\lambda\rightarrow v $ weakly in ${\rm X}$ and
strongly in ${\rm L}^q(\Omega)$ for any
$q\in\big[1,\frac{2N}{N-2}\big).$ By \eqref{e1} we deduce that $v
$ solves $-\Delta v=\lambda_j v$ in $\Omega,$ $v=0$ on
$\partial\Omega $ and, so, $v\in K_j.$

We claim that $v\not=0.$ Indeed, if $v=0$ by \eqref{e1} we get
$$
\int\limits_\Omega|\nabla v_\lambda|^2\,
dx=\lambda\int\limits_\Omega v_\lambda^2\,dx+ \|u_\lambda\|_{\rm
X}^{p-1}\int\limits_\Omega |v_\lambda|^{p+1}\,dx
$$
and passing to the limit as $\lambda$ goes to $\lambda_j$ we deduce
that $\|v_\lambda\|$ goes to zero, since $\|v_\lambda\|_{p+1}$ is
bounded. Moreover by \eqref{e1} and Lemma \ref{adn2} we get
$$
\|v_\lambda\|_s\le c\left(\lambda\|v_\lambda\|_\frac{Ns}{N+2s}
+\|u_\lambda\|_{\rm X}^{p-1}
 \|v_\lambda\|_\frac{pNs}{N+2s}^p\right),
$$
and by interpolation (since $1<\frac{Ns}{N+2s}<s$)
$\|v_\lambda\|_s\rightarrow 0$ as $\lambda$ goes to $\lambda_j$.

Finally a contradiction arises since $\|v_\lambda\|_{\rm X}=1.$

Now it is easy to check that there exist $b_\lambda^1,\dots,
b_\lambda^k\in\R$ and $\psi_\lambda\in K_j^\perp$ such that
$v_\lambda=\sum\limits_{i=1}^kb^i_\lambda e_i+\psi_\lambda$ and, as
$\lambda$ goes to $\lambda_j$, $\|\psi_\lambda\|_{\rm X}\rightarrow
0$ (using Lemma \ref{adn2}) and $b_\lambda^i\rightarrow b^i,$
$i=1,\dots,k$ (see also \cite{cocv} for analogous properties in
presence of more general nonlinearities).

We want to prove that there exists $\Lambda$ such that
\begin{equation}\label{e2}
\frac{\|u_\lambda\|_{\rm X}^{p-1}}{\lambda_j-\lambda}\rightarrow
\Lambda>0\ \hbox{as $\lambda$ goes to $\lambda_j$.}
\end{equation}
Multiplying \eqref{e1} by $e_i,$ we deduce that
\begin{equation}\label{e3}
b^i_\lambda(\lambda_j-\lambda)=\|u_\lambda\|_{\rm
X}^{p-1}\int\limits_\Omega  |v_\lambda|^{p-1}v_\lambda
e_i\,dx,\quad i=1,\dots,k.
\end{equation}
We recall that,  as $\lambda$ goes to $\lambda_j$,
$b_\lambda^i\rightarrow b^i $ and $b^i\not=0$ for some $i$, since
$v\not=0$. We also point out that, as $\lambda$ goes to
$\lambda_j$, $\int\limits_\Omega |v_\lambda|^{p-1}v_\lambda e_i
\rightarrow \int\limits_\Omega |v|^{p-1}v e_i$ (since
$v_\la\rightarrow v$ strongly in $X$), and $\int\limits_\Omega
|v|^{p-1}v e_i\not=0$ for some $i$ (since $\int\limits_\Omega
|v|^{p+1}\not=0$). Then by \eqref{e3} we deduce that, as $\lambda$
goes to $\lambda_j$, $\frac{\|u_\lambda\|_{\rm
X}^{p-1}}{\lambda_j-\lambda}$ goes to $\Lambda\in\R, $
$\Lambda\not=0 $ and also that
\begin{equation}\label{e4}
b^i=\Lambda\int\limits_\Omega |v|^{p-1}v e_i\,dx=
\Lambda\int\limits_\Omega\left|\sum\limits_{h=1}^k
b^he_h\right|^{p-1} \left(\sum\limits_{h=1}^k b^he_h\right)
e_i\,dx.
\end{equation}
By \eqref{e4} we easily deduce  that $\Lambda>0$ and \eqref{e2} is
proved.

Finally, we can write $u_\lambda=(\lambda_j-\lambda)^\frac{1}{p-1}
\left(\sum\limits_{i=1}^ka_\lambda^ie_i+\phi_\lambda\right),$
where $a^i_\lambda=  \|u_\lambda\|_{\rm X}
(\lambda_j-\lambda)^{-\frac{1}{p-1}}b^i_\lambda$ and
$\phi_\lambda= \|u_\lambda\|_{\rm X}
(\lambda_j-\lambda)^{-\frac{1}{p-1}}\psi_\lambda$. Moreover, as
$\lambda$ goes to $\lambda_j$, $a_\lambda$ goes to
$a=\Lambda^\frac{1}{p-1} b $ and by \eqref{e4} we deduce that $a$
is a critical point of $J_{\lambda_j}.$ That proves our
claim.\end{proof}

Secondly, we prove that any non degenerate critical point $a$ of
$J_{\lambda_j}$ generates a unique solution $u_\lambda$
bifurcating from the eigenvalue $\lambda_j$ which satisfies
\eqref{ul}.

\begin{prop}\label{dege}
Suppose that $a$ is a non degenerate nontrivial critical point of
the function $J_{\lambda_j}$ defined in \eqref{J}. Then there
exists $\delta=\delta(\la_j)>0$ such that for any
$\lambda\in(\lambda_j-\delta,\lambda_j)$ problem \eqref{p2} with
$\epsilon=\la_j-\la$ has a unique solution $u_\lambda$ such that
$u_\lambda=a_\lambda e+\phi_\lambda,$ where $a_\lambda \to a$ in
$\R^k$, $\langle \phi_\lambda, e_i\rangle=0$ for any
$i=1,\ldots,k$ and $\|\phi_\lambda\|_{\rm X}\to 0$ as
$\lambda\to\lambda_j$.
\end{prop}

\begin{proof}

$\underline{\hbox{ Step 1. Existence}}$

Since $a$ is a non degenerate critical point of  $J_{\lambda_j}$,
by Lemma \ref{jj} we deduce  that there exists $\delta>0$ such
that for any $\lambda\in(\lambda_j-\delta,\lambda_j)$  the
function $J_\lambda$ has a critical point $a_\lambda$ such that
$a_\lambda$ goes to $a$ as as $\lambda$ goes to $\lambda_j$. Then
by Lemma \ref{jt} we deduce that the function
$u_\lambda=a_\lambda e+\phi_\lambda(a_\lambda) $ is a solution  to
problem \eqref{p2}, with $\langle\phi_\lambda, e_i\rangle=0$ for
any $i=1,\ldots,k$ and $\|\phi_\lambda\|_{\rm X}\to 0$ as
$\lambda\to\lambda_j$.

$\underline{\hbox{ Step 2. Uniqueness}}$

Let $u_\lambda $ and $v_\lambda$ be two solutions of \eqref{p2}
such that $u_\lambda=a_\lambda e+\phi_\lambda$ and
$v_\lambda=b_\lambda e+\psi_\lambda$, where
$\phi_\lambda,\psi_\lambda\in K_j^\perp,$ $a_\lambda,b_\lambda$ go
to $a$ and $\|\phi_\lambda\|_{\rm X},\|\psi_\lambda\|_{\rm X} $ go
to zero as $\lambda$ goes to $\lambda_j$.

Assume by contradiction that $u_ \lambda  \neq v_\lambda$ and
consider the function
\[
z_\lambda:=\dfrac{u_\lambda-v_\lambda}{\|u_\lambda-v_\lambda\|}.
\]
It is clear that $z_\lambda$ satisfies the problem
\begin{equation}\label{pbze}
\begin{cases}
-\Delta z_\lambda=\lambda  z_\lambda+(\lambda_j-\la)
\dfrac{f(u_\lambda)-f(v_\lambda)}{\|u_\lambda-v_\lambda\|}
& \mbox{ in }\Omega\\
z_\lambda=0 & \mbox{ on }\partial \Omega.
\end{cases}
\end{equation}

We point out that  by the Mean Value Theorem there exists
$\vartheta\in(0,1)$ such that
\begin{equation}\label{mvt}
\dfrac{f(u_\lambda)-f(v_\lambda)}{\|u_\lambda-v_\lambda\|}=f'
\left(u_\lambda+\vartheta(u_\lambda-v_\lambda)\right)z_\lambda
.\end{equation}
We also remark that
$f'\left(u_\lambda+\vartheta(u_\lambda-v_\lambda)\right)$
converges to $f'(ae)$ strongly in ${\rm L^{N/2}}(\Omega)$  as
$\lambda$ goes to $\lambda_j$.

Up to a subsequence, we can assume that $z_\lambda\rightarrow z$
weakly in ${\rm H}^1_0(\Omega)$ and strongly in ${\rm
L}^q(\Omega)$ for any $1<q<\frac{2N}{N-2}.$ Moreover, by
\eqref{pbze} we deduce that there exists
$\alpha=(\alpha_1,\ldots,\alpha_k)\in \R^k$ such that $z=\alpha
e=\sum_{i=1}^k \alpha_ie_i$.

Now, multiplying \eqref{pbze} by $e_i$, $i=1,\ldots,k$, and  using
(\ref{mvt}), we deduce
$$
\int\limits_\Omega z_\lambda e_i\, dx=\int\limits_\Omega
f'\left(u_\lambda+\vartheta(u_\lambda-v_\lambda)\right)z_\lambda
e_i\,dx
$$
and passing to the limit, as $\lambda$ goes to $\lambda_j$, we get
$\alpha_i=\int\limits_\Omega f'(ae)z e_i\,dx$ for any
$i=1,\ldots,k.$ Therefore $\alpha$ is a solution of the linear
system ${\cal H} { J_{\lambda_j}}(a)\alpha=0, $ where $\mathcal{H}
{J_{\lambda_j}}(a)$ denotes the Hessian matrix of $\tilde{J}$ at
$a$. Since $a$ is a non degenerate critical point of
$J_{\lambda_j}$, we deduce that $\alpha=0$, namely $z=0.$

On the other hand, multiplying \eqref{pbze} by $z_\lambda$, and
using (\ref{mvt}), we deduce
$$
\int\limits_\Omega |\nabla z_\lambda|^2\,dx=
\lambda\int\limits_\Omega  z_\lambda ^2\,dx
+(\lambda_j-\lambda)\int\limits_\Omega
f'\left(u_\lambda+\vartheta(u_\lambda-v_\lambda)\right)z_\lambda
^2\,dx
$$
and passing to the limit, as $\lambda$ goes to $\lambda_j$, we get
$\|z_\lambda\|\rightarrow 0$. Finally, a contradiction arises
since $\|z_\lambda\|=1.$
\end{proof}

Finally, we compute the Morse index of the solution $u_\lambda $
generated by a critical point $a$ of $J_{\lambda_j}$ (see
\eqref{J}) in terms of the Morse index of $a$.

We recall that the Morse index of a solution $u$ of problem
(\ref{p1}) is the number of negative eigenvalues $\mu$ of the
linear problem
$$
v-i^*\left[\lambda v+ f'(u)v\right]=\mu v,\ v\in{\rm H}^1_0(\Omega)
$$
or equivalently
$$
\begin{cases}
-(1-\mu)\Delta v=\lambda v+ f'(u)v &  \hbox{ in $\Omega,$}\\
v=0 & \hbox{ on $\partial\Omega.$}\end{cases}
$$

We point out that the  function $u$, which solves  problem
(\ref{p1}), and the function $v=\eps^{-\frac{1}{p-1}}u$, which
solves problem (\ref{p2}), have the same Morse index.

\begin{prop}\label{morse}
Let $u_\lambda=\sum\limits_{i=1}^ka^i_\lambda e_i+\phi_\lambda$ be
a solution to \eqref{p2} such that
$\lim\limits_{\la\to\la_j}\|\phi_\la\|_{\rm X}=0$,
$\lim\limits_{\lambda\rightarrow \lambda_j}a_\lambda^i=a^i$ and
$(a^1,\dots,a^k)$ is a non trivial critical point of
$J_{\lambda_j}
 $ (see \eqref{J}). If the Morse index of $a$ is $m,$  then the
Morse index of $u_\lambda$ is at least $m+j-1.$ Moreover if $a$ is
also non degenerate, then the solution $u_\lambda$ is non
degenerate and its Morse index is exactly $m+j-1.$
\end{prop}

\begin{proof} We denote by
$\mu_\lambda^1<\mu_\lambda^2\le\dots\le\mu_\lambda^i\le\dots$ the
sequence of the eigenvalues, counted with their multiplicities, of
the linear problem
\begin{equation}\label{linea}
\begin{cases} -(1-\mu)\Delta v=\lambda v+(\lambda_j-\lambda)f'(u_\la)v
& \hbox{ in $\Omega,$}\\ v=0 & \hbox{ on
$\partial\Omega.$}\end{cases}
\end{equation}

We also denote by $v_\lambda^i\in{\rm H}^1_0(\Omega),$ with
$\|v_\lambda^i\|_2=1,$ the eigenfunction associated to the
eigenvalue $\mu_\lambda^i.$

It is clear that, as $\lambda$ goes to $\lambda_j,$ eigenvalues
and eigenfunctions of  (\ref{linea}) converge to eigenvalues and
eigenfunctions of the linear problem
$$
\begin{cases}
-(1-\mu)\Delta v=\lambda_j v &  \hbox{ in $\Omega,$}\\ v=0 &
\hbox{ on $\partial\Omega,$}\end{cases}
 $$
whose set of eigenvalues is
$$
\left\{1-\frac{\lambda_j}{\lambda_1},\dots,1-\frac{\lambda_j}{\lambda_{j-1}},
{\mathop{ \underbrace{0,\dots,0}
}\limits_k},1-\frac{\lambda_j}{\lambda_{j+k }},\dots\right\}.
$$

Therefore, if $\lambda$ is close enough to $\lambda_j$, we can claim
that $\mu_\lambda^1,\dots,\mu_\lambda^{j-1}$ are negative and they
are close to $1-\frac{\lambda_j}{\lambda_1},\dots,1-
\frac{\lambda_j}{\lambda_{j-1}},$ respectively, and that
$\mu_\lambda ^{j+k }$ is positive and close to
$1-\frac{\lambda_j}{\lambda_{j+k }}$.

Therefore, it remains to understand what happens to the $k$
eigenvalues $\mu_\lambda^j,\dots,\mu_\lambda^{j+k-1}$, which go to
zero as $\lambda$ goes to $\lambda_j.$

We claim that
 \begin{equation}\label{eige2}
 \begin{cases}
\lim\limits_{\lambda\rightarrow\lambda_j }
\frac{\mu^{j+l-1}_\lambda}{\lambda-\lambda_j}\lambda_j=\Lambda^l,\quad
\hbox{where } \Lambda^1\le\dots\le\Lambda^k,\ l=1,\dots, k, \mbox{ and}\\
\hbox{ $\Lambda^l$ are the eigenvalues of the Hessian matrix
${\cal H}J_{\lambda_j}(a).$}
\end{cases}
\end{equation}

For any $l=1,\dots,k$ we denote by $v^l_\lambda$ an eigenfunction
associated to $\mu^{j+l-1}_\lambda,$ with $\| v^l_\lambda\|_2=1,$
i.e.
\begin{equation}\label{linj}
\begin{cases}
-(1-\mu^{j+l-1}_\lambda)\Delta v^l_\lambda=\lambda
v^l_\lambda+(\lambda_j-\lambda)f'(u_\lambda)v^l_\lambda &  \hbox{
in $\Omega,$}\\ v^l_\lambda=0 & \hbox{ on
$\partial\Omega.$}\end{cases}
\end{equation}

Then we can write
\begin{equation}\label{eige1}
\begin{cases}
v^l_\lambda=\sum\limits_{i=1}^k
b_\lambda^{l,i}e_i+\psi^l_\lambda,\  b_\lambda^{l,i}\in\R,\\
\langle\psi^l_\lambda,e_i\rangle=0,\ i=1,\dots,k,\ \langle
v^l_\lambda,v^s_\lambda\rangle=0\ \hbox{if}\
l\not=s,  \\
\sum\limits_{i=1}^k
\left(b_\lambda^{l,i}\right)^2+\|\psi_\la^l\|^2_2=1.
\end{cases}
\end{equation}

Now, up to a subsequence, we can assume that for any $l=1,\dots,k$
and $i=1,\dots,k$, $\psi_\lambda^l\rightarrow\psi^l$ and $
b_\lambda^{l,i}\rightarrow  b ^{l,i}$ as $\lambda$ goes to
$\lambda_j.$ Then $v^l_\lambda\rightarrow v^l:=\sum\limits_{i=1}^k
b ^{l,i}e_i+\psi^l $ as $\lambda$ goes to $\lambda_j.$ We point
out that the convergence in ${\rm H}^1_0(\Omega)$ is strong, since
$v_\lambda^l$ solves equation (\ref{linj}) and
$\mu_\lambda^{j+l-1}$ does not go to 1 as $\lambda$ goes to
$\lambda_j.$

First of all we claim that $\psi^l=0$ for any $l=1,\dots,k$. In
fact by (\ref{linj}) we deduce that  for any $l=1,\dots,k$ and for
all $v\in{\rm H}^1_0(\Omega)$ it holds
$$
(1-\mu^{j+l-1}_\lambda)\int\limits_\Omega \nabla v^l_\lambda\nabla
v\,dx=\lambda \int\limits_\Omega v^l_\lambda
v\,dx+(\lambda_j-\lambda)\int\limits_\Omega
f'(u_\lambda)v^l_\lambda v\,dx,
$$
and passing to the limit as $\lambda$ goes to $\lambda_j $
$$
\int\limits_\Omega \nabla v^l \nabla v\,dx=\lambda_j
\int\limits_\Omega v^l_\lambda  v\,dx \quad\forall\ v\in{\rm
H}^1_0(\Omega),
$$
that is $v^l$ is an eigenfunction associated to the eigenvalue
$\lambda_j,$ so that $\psi^l=0.$ Therefore by (\ref{eige1}) we
deduce that
\begin{equation}\label{eige11}
\sum\limits_{i=1}^k  b ^{l,i}b^{s,i}=0 \ \hbox{if $l\not=s$ and }
\sum\limits_{i=1}^k  \left(b ^{l,i}\right)^2=1 \quad
\forall\,l=1,\ldots,k.
\end{equation}

Now, if we multiply (\ref{linj}) by $e_i$, we get for any
$i=1,\dots,k$ and for any  $l=1,\dots,k$,
$$
(1-\mu^{j+l-1}_\lambda)\lambda_jb^{l,i}_\lambda =\lambda
b^{l,i}_\lambda+ (\lambda_j-\lambda)\int\limits_\Omega
f'(u_\lambda)v^l_\lambda e_i\,dx,
$$
that is
\begin{equation}\label{eige3}
b^{l,i}_\lambda
\left(1-\frac{\mu^{j+l-1}_\lambda}{\lambda_j-\lambda}\lambda_j\right)=
\int\limits_\Omega f'(u_\lambda) \left(\sum\limits_{i=1}^k
b_\lambda^{l,i}e_i+\psi^l_\lambda\right)e_i\,dx.
\end{equation}
Now, since as $\lambda$ goes to $\lambda_j$
\begin{eqnarray*}
\int\limits_\Omega f'\left(\sum\limits_{i=1}^k a_\lambda^{
i}e_i+\phi _\lambda\right) \left(\sum\limits_{i=1}^k
b_\lambda^{l,i}e_i+\psi^l_\lambda\right)e_i\,dx\\
\rightarrow
\int\limits_\Omega f'\left(\sum\limits_{i=1}^k a ^{ i}e_i \right)
\left(\sum\limits_{i=1}^k b ^{l,i}e_i \right)e_i\,dx,
\end{eqnarray*}
by (\ref{eige3}) and by (\ref{eige11}) we deduce that for any
$l=1,\dots,k$, when $\lambda$ goes to $\lambda_j$,
$\frac{\mu^{j+l-1}_\lambda}{\lambda_j-\lambda}\lambda_j$ converges
to an eigenvalue $\Lambda^l$ of the matrix ${\cal H}J(a)$ and also
that $b^l_\lambda e$ converges to the associated eigenfunction
$b^le,$ since
$$
b^{l,i}  \left(1-\Lambda^l\right)= \int\limits_\Omega
f'\left(\sum\limits_{i=1}^k a ^{ i}e_i \right)
\left(\sum\limits_{i=1}^k b ^{l,i}e_i \right)e_i\,dx \quad
\forall\,i=1,\ldots,k,
$$
i.e. $HJ(a)(b^le)=\Lambda^l(b^le)$.

Since $a$ has Morse index $m$, there are $m$ eigenvalues
$\Lambda^l$ which are negative, so that at least $m$ eigenvalues
$\mu_\la^l$ are negative as well, provided $\la$ is close to
$\la_j$.

Finally, if $a$ is non degenerate, all the $\Lambda^l$'s are
different from 0, so that, if $\la$ is near $\la_j$, there are
exactly $m$ negative eigenvalues $\mu_\la^l$, as claimed.
\end{proof}

\begin{proof}[Proof of Theorem \ref{main2}].

The claim follows by Proposition \ref{dege} and Proposition
\ref{morse}.
\end{proof}

\begin{proof}[Proof of Theorem \ref{main3}]

The claim follows by Proposition \ref{nece} and Theorem
\ref{main2}.
\end{proof}

\sezione{ Examples}

In the first example we compute the Morse index of the solutions
bifurcating from a simple eigenvalue.

\begin{ex}\label{semplice}
If $\lambda_j$ is a simple eigenvalue, then there exists $\delta
>0$ such that for any $\lambda\in (\lambda_j-\delta,\lambda_j)$
the problem \eqref{p2} has exactly one pair of solutions
$(u_\lambda,-u_\lambda)$ bifurcating from $\lambda_j$, whose Morse
index is $j.$
\end{ex}

\begin{proof} In this particular case the function $J_{\lambda_j}$
reduces to $J_{\lambda_j}(a)=\frac{1}{2}
a^2-\frac{1}{p+1}|a|^{p+1}\int\limits_\Omega|e_1|^{p+1}$. It is
easy to prove that $J_{\lambda_j}$ has exactly two (non trivial)
critical points,
$a_0:=\left(\int\limits_\Omega|e_1|^{p+1}\right)^{-\frac{1}{p-1}}$
and $-a_0,$ which are non degenerate and have Morse index 1. The
claim follows by Theorem \ref{main2} and Theorem \ref{main3}.

\end{proof}

In the second example we study solutions bifurcating from any
multiple eigenvalues when $\Omega$ is a rectangle in $\R^2.$.

\begin{ex}\label{rect}
Let $\Omega=(0,L)\times(0,M)$ be a rectangle in $\R^2.$ Let
$\lambda_j $ be an eigenvalue for $\Omega,$ which has multiplicity
$k$. Then there exists $\delta >0$ such that for any $\lambda\in
(\lambda_j-\delta,\lambda_j)$ problem \eqref{p2} with $p=3$ has
exactly $\frac{3^k-1}{2}$ pairs of solutions
$(u_\lambda,-u_\lambda)$ bifurcating from $\lambda_j$. In
particular, if $k=2$ problem \eqref{p2} has two pairs of solutions
with Morse index $j+1$ and two pairs of solutions with Morse index
$j $ and if $k=3$ problem \eqref{p2} has three pairs of solutions
with Morse index $j+2$, six pairs of solutions with Morse index
$j+1$ and four pairs of solutions with Morse index $j$.
\end{ex}

\begin{proof} We know that $\lambda_j=\pi^2\left(\frac{n_i^2}{L^2}+
\frac{m_i^2}{M^2}\right) $, $i=1,\dots,k$, where $n_i,m_i\in
\mathbb{N}$ with $n_i\not=n_l$ and $m_i\not=m_l $ if $i\not= l.$
The eigenspace associated to $\lambda_j$ is spanned by the
functions
$$
e_i(x,y)=\sqrt{\frac{4}{LM}}\sin \frac{n_i\pi}{L}x\sin
\frac{m_i\pi}{M}y,\ i=1,\dots,k.
$$

Taking in account that
$$
\int\limits_\Omega e_i^3e_l\,dx=\int\limits_\Omega
e_i^2e_he_l\,dx=0\ \hbox{if}\ i,h\not= l
$$
and also that
$$
\alpha:=\int\limits_\Omega e_i^4\,dx=\frac{9}{64}LM\ \hbox{ and }\
\beta:=\int\limits_\Omega e_i^2e_l^2\,dx=\frac{1}{16}LM,
$$
the function  $J_{\lambda_j}:\R^k\rightarrow\R$ introduced in
(\ref{J}) with $p=3$ reduces to
$$
J_{\lambda_j}(a )=\frac{1}{2}\sum\limits_{i=1}^ka_i^2
-\frac{1}{4}\alpha\sum\limits_{i=1}^ka_i^4 -\frac{3}{4}\beta
\sum\limits_{i,l=1\atop i\not=l}^ka_i^2a_l^2.
$$

We can compute
\begin{equation}\label{grad}
\frac{\partial J_{\lambda_j}}{\partial a_i}(a  )=a_i-\alpha
a_i^3-3\beta a_i \sum\limits_{ l=1\atop l\not=i}^ka_l^2 ,\
i=1,\dots,k \end{equation} and also the Hessian matrix ${\cal
H}J_{\lambda_j}(a )$
\begin{equation}\label{matrice}
\begin{pmatrix}
1-3\alpha a_1^2-3\beta\sum\limits_{ l=1\atop l\not=1}^ka_l^2&
-6\beta a_1 a_2 & \dots& -6\beta a_1 a_k \\
-6\beta a_1 a_2 & 1-3\alpha a_2^2-3\beta\sum\limits_{ l=1\atop l\not=2}^ka_l^2
& \dots & -6\beta a_2 a_k\\
\vdots & \vdots & \ddots & \vdots\\
-6\beta a_1 a_k & -6\beta a_2 a_k& \dots&  1-3\alpha
a_k^2-3\beta\sum\limits_{ l=1\atop l\not=k}^ka_l^2
\end{pmatrix}.
\end{equation}

Let us consider the case $k=2.$ It is easy to check that
(\ref{grad})reduces to
\begin{eqnarray*}
& &\frac{\partial J_{\lambda_j}}{\partial a_1}(a_1,a_2
)=a_1-\alpha
a_1^3-3\beta a_1  a_2^2 \\
& &\frac{\partial J_{\lambda_j}}{\partial a_2}(a_1,a_2
)=a_2-\alpha
a_2^3-3\beta a_2  a_1^2 \\
\end{eqnarray*}
and also that (\ref{matrice}) reduces to
$$
{\cal H}
J_{\lambda_j}(a_1,a_2)=\begin{pmatrix}
1-3\alpha a_1^2-3\beta a_2^2&-6\beta a_1a_2 \\
-6\beta a_1a_2  & 1-3\alpha a_2^2-3\beta a_1^2
\end{pmatrix}.
$$
Let $\gamma_1:=\alpha^{-1/2} $ and
$\gamma_2:=(\alpha+3\beta)^{-1/2}.$ Therefore  $J_{\lambda_j}$ has
exactly the following (non trivial)  critical points: $(0,\pm
\gamma_1),$ $(\pm\gamma_1,0 ),$ which have Morse index 2 and
$(\gamma_2,\pm\gamma_2)$ and $(-\gamma_2,\pm\gamma_2),$ which have
Morse index 1. Therefore the claim follows by Theorem \ref{main2}.

Let us consider the case $k=3.$ It is easy to check that
(\ref{grad})reduces to
\begin{eqnarray*}
& &\frac{\partial J_{\lambda_j}}{\partial a_1}(a_1,a_2,a_3)=
a_1-\alpha a_1^3-3\beta a_1(a_2^2+a_3^2)\\
& &\frac{\partial J_{\lambda_j}}{\partial a_2}(a_1,a_2,a_3)=
a_2-\alpha a_2^3-3\beta a_1(a_1^2+a_3^2)\\
& &\frac{\partial J_{\lambda_j}}{\partial
a_3}(a_1,a_2,a_3)=a_3-\alpha a_3^3-3\beta a_1(a_1^2+a_2^2)
\end{eqnarray*} and also that the Hessian matrix ${\cal H}
J_{\lambda_j}(a_1,a_2,a_3)$
(\ref{matrice})reduces to
$$
\begin{pmatrix}
1-3\alpha a_1^2-3\beta(a_2^2+a_3^2) & -6\beta a_1 a_2& -6\beta a_1 a_3 \\
-6\beta a_1 a_2 & 1-3\alpha a_2^2-3\beta(a_1^2+a_3^2)&-6\beta a_2 a_3\\
-6\beta a_1 a_3 &-6\beta a_2 a_3&  1-3\alpha
a_3^2-3\beta(a_1^2+a_2^2)
\end{pmatrix}.
$$
Let $\gamma_1:=\alpha^{-1/2},$ $\gamma_2:=(\alpha+3\beta)^{-1/2}$
and $\gamma_3:=(\alpha+6\beta)^{-1/2}.$ Therefore $J_{\lambda_j}$
has exactly the following (non trivial) pairs of critical points:
\begin{eqnarray*}
& &A_1:=(0,0, \gamma_1),\ A_2:=(0,  \gamma_1,0 ), \ A_3:=( \gamma_1,0,0)\\
& &A_4:=(0,\gamma_2, \gamma_2),\ A_5:=(0,\gamma_2, -\gamma_2),\
A_6:=( \gamma_2, 0,\gamma_2),\\
& &  A_7:=( \gamma_2,0, -\gamma_2),\ A_8:=( \gamma_2,
\gamma_2,0),\
A_9:= ( \gamma_2,  -\gamma_2,0), \\
& &A_{10}:=(\gamma_3,\gamma_3,\gamma_3),\
A_{11}:=(\gamma_3,\gamma_3,-\gamma_3),\\
& & A_{12}:=(\gamma_3,-\gamma_3,\gamma_3),\
A_{13}:=(-\gamma_3,\gamma_3,\gamma_3)
\end{eqnarray*}
and $-A_i$ for $i=1,\dots,13.$ A simple computation shows that
$A_i$ and $-A_i$ have Morse index 3 if $i=1,2,3,$ that they have
Morse index 2 if $i=4,\dots,9$ and that they have Morse index 1 if
$i=10,\dots,13.$ Therefore the  claim follows by Theorem
\ref{main2}.

In the general case, the function $J_{\lambda_j}$ has $3^k-1$
critical points of the form
$(0,\dots,0,\sotto\limits_k,0,\dots,0)$ which are non degenerate,
where $\gamma_i:=\left[\alpha+3(i-1)\beta\right]^{-1/2}$ for
$i=1,\dots,k.$ Let us compute the Hessian matrix ${\cal H}
J_{\lambda_j}(a)$ (\ref{matrice}) at the point $a=(0,\ldots,0,
\sotto\limits_k,0,\dots,0).$ We have
$$
\frac{1}{\left[\alpha+3(i-1)\beta\right]^{k}}\begin{pmatrix}
A_i& 0\\
0 &B_i
\end{pmatrix},
$$
where the $i\times i$ matrix $A_i$ is given by
$$
A_i:=\begin{pmatrix}
-2\alpha & -6\beta  & \dots &-6\beta  \\
-6\beta   & - 2\alpha    & \dots     &-6\beta  \\
\vdots&\vdots&\ddots&\vdots\\
 -6\beta   &-6\beta  &  \dots &
-2\alpha
\end{pmatrix}=(ab)^i\begin{pmatrix}
-\frac{9}{32} & -\frac{3}{8} & \dots &-\frac{3}{8} \\
-\frac{3}{8}   & - \frac{9}{32}    & \dots     &-\frac{3}{8}  \\
\vdots&\vdots&\ddots&\vdots\\
 -\frac{3}{8}   &-\frac{3}{8} &  \dots &
-\frac{9}{32}
\end{pmatrix}
$$
and the $(k-i)\times (k-i)$ matrix $B_i$ is given by
$$
B_i:=\begin{pmatrix}
 \alpha  -3\beta  &0 & \dots &0  \\
0   & \alpha  -3\beta    & \dots     &0 \\
\vdots&\vdots&\ddots&\vdots\\
0  &0  &  \dots & \alpha  -3\beta
\end{pmatrix}=\begin{pmatrix}
 -\frac{3}{64}   &0 & \dots &0  \\
0   & -\frac{3}{64}    & \dots     &0 \\
\vdots&\vdots&\ddots&\vdots\\
0  &0  &  \dots & -\frac{3}{64}
\end{pmatrix}.
$$
\end{proof}

In the third  example we study solutions bifurcating from the second
eigenvalue  when $\Omega$ is a cube in $\R^3.$

\begin{ex}\label{cubo}
Let $\Omega=(0,\pi)\times(0,\pi)\times(0,\pi)$ be the cube in
$\R^3.$ In this case $\lambda_2=6$ is the second eigenvalue for
$\Omega$ and it has multiplicity three. Then there exists $\delta
>0$ such that for any $\lambda\in (\lambda_2-\delta,\lambda_2)$
problem \eqref{p2} with $p=3$ has exactly $13$ pairs of  solutions
$(u_\lambda,-u_\lambda)$ bifurcating from $\lambda_2$. Moreover
three pairs of solutions have Morse index $3$, six pairs of
solutions have Morse index $2$ and four pairs of solutions have
Morse index $1$.
\end{ex}

\begin{proof}

The eigenspace associated to $\lambda_2$ is spanned by the functions
\begin{eqnarray*}
& &e_1(x,y,z)=\sqrt{\frac{8}{\pi^3}}\sin x\sin y\sin 2z,\\
& & e_2(x,y,z)=\sqrt{\frac{8}{\pi^3}}\sin x\sin2 y\sin  z,\\
& &e_1(x,y,z)=\sqrt{\frac{8}{\pi^3}}\sin 2x\sin y\sin  z.
\end{eqnarray*}

Taking in account that
$$
\int\limits_\Omega e_1^3 e_2dx=\int\limits_\Omega e_1^2 e_2e_3dx=0
$$
and also that
$$
\alpha:=\int\limits_\Omega e_1^4dx=\frac{27}{8\pi^3}\ \hbox{and}\
\beta:=\int\limits_\Omega e_1^2e_2^2dx=\frac{3}{2\pi^3},
$$
the function $J_{\lambda_j}$ with $p=3$ reduces to
$$
J_{\lambda_j}(a_1,a_2,a_3)=\frac{1}{2}(a_1^2+a_2^2+a_3^2)-\frac{1}{4}
\alpha(a_1^4+a_2^4+a_3^4)-\frac{3}{2}\beta
(a_1^2a_2^2+a_1^2a_3^2+a_2^2a_3^2).
$$

Therefore the  claim follows by Theorem \ref{main2}, arguing
exactly as in the previous example in the case $k=3.$ \end{proof}



\end{document}